%% file: 2005-51.tex
\documentclass{gtart_h}  

\input gtoutput

\lognumber{552}
\received{14 January 2005}
\volumenumber{9}\papernumber{51}\volumeyear{2005}
\pagenumbers{2227}{2259}   
\published{3 December 2005}
\accepted{26 November 2005}
\proposed{Cameron Gordon}
\seconded{David Gabai, Joan Birman}

\usepackage{amssymb,amsmath,amscd,graphicx,color}
\def\SetFigFont#1#2#3#4#5{\small}

\def\figref#1{\hyperlink{#1anchor}{Figure~\ref*{#1}}}
\def\anchor#1{\noindent\hypertarget{#1anchor}{\smash{$\phantom{99}$}}}

\let\tilde\widetilde

\def\a{\alpha}
\def\A{\mathcal{A}}
\def\b{\beta}
\def\g{\gamma}
\def\G{\Gamma}
\def\E{\mathcal{E}}
\def\D{\Delta}
\def\H{\mathbb{H}}

\def\l{\lambda}

\def\R{\mathbb{R}}
\def\s{\sigma}
\def\S{\Sigma}

\def\d{\partial}
\def\cross{\times}

\def\e{\epsilon}
\def\half{\frac{1}{2}}

\def\vol{\text{vol}}

\newtheorem{theorem}{Theorem}[section]
\newtheorem{lemma}[theorem]{Lemma}
\newtheorem{corollary}[theorem]{Corollary}

\theoremstyle{definition}

\newtheorem{definition}[theorem]{Definition}

\begin{document}

\title{Heegaard gradient and virtual fibers}
\author{Joseph Maher}
\address{Mathematics 253-37, California Institute of
Technology\\Pasadena, CA 91125, USA}

\email{maher@its.caltech.edu}

\begin{abstract}
We show that if a closed hyperbolic $3$--manifold has infinitely many
finite covers of bounded Heegaard genus, then it is virtually
fibered. This generalizes a theorem of Lackenby, removing restrictions
needed about the regularity of the covers. Furthermore, we can replace
the assumption that the covers have bounded Heegaard genus with the
weaker hypotheses that the Heegaard splittings for the covers have
Heegaard gradient zero, and also bounded width, in the sense of
Scharlemann--Thompson thin position for Heegaard splittings.
\end{abstract}

\asciiabstract{%
We show that if a closed hyperbolic 3-manifold has infinitely many
finite covers of bounded Heegaard genus, then it is virtually
fibered. This generalizes a theorem of Lackenby, removing restrictions
needed about the regularity of the covers. Furthermore, we can replace
the assumption that the covers have bounded Heegaard genus with the
weaker hypotheses that the Heegaard splittings for the covers have
Heegaard gradient zero, and also bounded width, in the sense of
Scharlemann-Thompson thin position for Heegaard splittings.}

\keywords{Heegaard splitting, virtual fiber, hyperbolic $3$--manifold}
\asciikeywords{Heegaard splitting, virtual fiber, hyperbolic 3-manifold}

\primaryclass{57M10}\secondaryclass{57M50}

\maketitle

\section{Introduction}

Inspired by Lubotzky's work on Property ($\tau$) \cite{lubotzky}, Lackenby \cite{lac1,lac2} investigated the behaviour of Heegaard splittings in finite covers of hyperbolic $3$--manifolds, and how this relates to other well-known conjectures about $3$--manifolds, such as Thurston's conjecture that a hyperbolic $3$--manifold has a finite cover which has positive Betti number, or Thurston's stronger conjecture that every hyperbolic $3$--manifold has a finite cover which is fibered. Lackenby was able to show that if a hyperbolic manifold $M$ has infinitely many finite covers $M_i$ of bounded Heegaard genus, and bounded irregularity, then there is a finite cover which is fibered. Here, bounded irregularity means that the index of the normalizer of $\pi_1 M_i$ in $\pi_1 M$ is bounded. In this paper we show how to remove the assumptions about the regularity of the covers, and we can also replace the assumption about bounded genus with a weaker hypothesis. To be precise, we prove the following theorem:

\begin{theorem} \label{main}
Let $X$ be a closed hyperbolic $3$--orbifold, with infinitely many finite manifold covers $M_i$ of degree $d_i$, with Heegaard splittings $H_i$, of Euler characteristic $\chi_i$. If the Heegaard gradient $\chi_i / d_i$ tends to zero, and the widths $c_+(H_i)$ of the Heegaard splittings are bounded, then all but finitely many $M_i$ contain an embedded surface which is a virtual fiber. Furthermore, there are only finitely many choices for the virtual fiber.
\end{theorem}

Here $c_+(H_i)$ is the width of the Scharlemann--Thompson thin position for the Heegaard splitting, which in particular, is bounded by the genus of the splitting. In fact, the proof only requires the width $c_+(H_i)$ to grow sufficiently slowly, but the growth bound we produce is logarithmic, and not particularly convenient to write down.

Note that if a collection of finite manifold covers $M_i$ of $X$ are not Haken, then the width of the splitting is the same as the Heegaard genus for each $M_i$, so if the widths are bounded, then the Heegaard gradient tends to zero. This gives the following corollary:

\begin{corollary}
Let $X$ be a closed hyperbolic $3$--orbifold with infinitely many finite covers $M_i$ for which $c_+(M_i)$ is bounded. Then $X$ is virtually Haken.
\end{corollary}

Ian Agol has also announced a proof of Theorem \ref{main}, and pointed out the following corollary.


\begin{corollary} \label{orbifold}
Let $X$ be a hyperbolic $3$--orbifold, and let $M$ be a manifold cover of $X$ of Heegaard genus at most $g$. There is a finite collection of manifold covers $C_1 \cup C_2 \cup C_3$ of $X$ such that one of the following occurs.

\begin{itemize}

\item $M \in C_1$

\item $M$ is a cyclic cover of $M' \in C_2$ dual to a fibration of $M'$ over $S^1$

\item $M$ is a subdihedral cover of $M' \in C_3$ dual to a fibration of $M'$ over $I^*$

\end{itemize}

\noindent Here $I^*$ is the $1$--orbifold formed by taking the quotient of $S^1$ by a reflection.
\end{corollary}

\begin{proof}
There are only finitely many covers of genus $g$ which do not contain an embedded virtual fiber, these make up the collection $C_1$. If $M$ contains an embedded virtual fiber $F$, then there is a common finite cover $\tilde M$ of $M$ and $X$, which is fibered. The pre-images of the components of $M - F$ in $\tilde M$ have fundamental group isomorphic to a closed surface group, so by a well-known result of Hempel, \cite[Theorem 10.6]{hempel}, the complementary pieces of $M - F$ are $I$--bundles over a surface. If there is a single connected complementary piece, then $M$ fibers over $S^1$, otherwise it fibers over $I^*$. As there are only finitely many choices for the virtual fiber, there is a finite collection of manifolds in each case for which all the other manifolds are cyclic or subdihedral covers.
\end{proof}

\subsection{Acknowledgements}

I would like to thank Ian Agol for informing me of Corollary \ref{orbifold}, and Nathan Dunfield for pointing out various mistakes in preliminary versions of this paper. I would also like to thank Marc Lackenby, Hyam Rubinstein, Jason Manning and Martin Scharlemann for helpful conversations. 

\subsection{Background}

Every closed $3$--manifold has a Heegaard splitting, which is a division of the manifold into two handlebodies. A handlebody is a compact $3$--manifold homeomorphic to a regular neighbourhood of a graph in $\R^3$. The graph is called a spine for the handlebody, and we say the genus of the handlebody is the genus of the surface forming its boundary.

A Heegaard splitting is a topological structure on a $3$--manifold, but it is possible to relate a Heegaard surface to a minimal surface, which is a geometric object in the manifold. A Heegaard splitting $S$ of a $3$--manifold $M$ gives rise to a sweepout, namely a degree one map $S \cross I \to M$, such that $S \cross 0$ and $S \cross 1$ are taken to the two spines of the handlebodies on either side of $S$. We will normally write $S_t$ for $S \cross t$. The area of the surfaces $S_t$ has some maximum value as $t$ varies. The minimax value is the smallest maximal area over all possible sweepouts of the manifold. If we take a sequence of maximal area sweepout surfaces whose area tends to the minimax value, then it is a fundamental result of Pitts and Rubinstein \cite{pr} (see also Colding and De Lellis \cite{cd}) that there is a subsequence which converges to a minimal surface in the manifold, possibly of lower genus. 

Recall that a minimal surface is a surface in a manifold with mean curvature zero. If the manifold is hyperbolic, then the intrinsic curvature of the surface is at most $-1$, so the Gauss--Bonnet formula gives an upper bound for the area of the minimal surface, and hence for any sweepout surface in a minimax sweepout. As some sweepout surface must divide the manifold into two parts of equal volume, this gives an immediate link between the Heegaard genus of the manifold and the Cheeger constant, which in turn is closely related to the first eigenvalue of the Laplacian. Lubotzky \cite{lubotzky} has shown that a closed hyperbolic $3$--manifold has Property ($\tau$) with respect to its congruence covers, which means that Cheeger constant of the covers is bounded away from zero. Lackenby \cite{lac1} uses this in his proof that the Heegaard genus of the congruence covers of a hyperbolic $3$--manifold grow linearly with their degree.

There is a bound on the diameter of a minimal surface in a hyperbolic $3$--manifold $M$ in terms of the injectivity radius of $M$, and the genus of the surface. As the injectivity radius of $M$ is a lower bound for the injectivity radius of any cover, this gives a diameter bound for minimal surfaces in covers of $M$ depending only on their genus. If the degree of the cover is very large with respect to the genus of the minimal surface, and the cover is fairly regular, then we can find many disjoint translates of the minimal surface in the cover. Lackenby is able to show the manifold is virtually fibered by showing that some of these translates must be parallel.

If we make no assumptions about the regularity of the covers, then we may not have any disjoint covering translations. However, if the volume of the manifold is large, and we have a sweepout of surfaces of bounded diameter, then we can find many disjoint sweepout surfaces in the manifold. The original sweepout surfaces need not have bounded diameter, but in Section \ref{section:generalized} we define a generalized sweepout, and show that we can find one with sweepout surfaces of bounded diameter. 

The diameter bound depends on the genus of the surface, so if infinitely many covers each have a Heegaard splitting of bounded genus, then they each have sweepouts of bounded genus. However, there is a useful technique called Scharl\-emann--Thompson thin position for Heegaard splittings \cite{st} which in some cases allows us to produce more general sweepouts with lower genus than the genus of the original Heegaard splitting. We now give a brief outline of this technique. We can think of a Heegaard splitting $H$ as a handle structure for the manifold, with one of the handlebodies constructed from a single zero-handle and a collection of one-handles, and the other handlebody constructed from a single three-handle and some two-handles. It may be possible to rearrange the order of the handles, ie, add some of the two-handles before all of the one-handles have been added. Such a re-ordering of the one- and two-handles determines a sequence of surfaces in the manifold, which divide the manifold into compression bodies consisting of handles of the same index which are consecutive in the ordering. Label the compression bodies corresponding to consecutive one-handles $A_1, \ldots A_{2n-1}$, and the compression bodies corresponding to consecutive two-handles $A_2, \ldots A_{2n}$, and let $F_i$ be the surface separating $A_{i-1}$ from  $A_i$. An even surface $F_{2i}$ is a Heegaard surface for the $3$--manifold consisting of the union of the adjacent pair of compression bodies, $A_{2i-1} \cup A_{2i}$. In \cite{st} Scharlemann and Thompson define the complexity of an ordering of the one-and two-handles to be the sequence of complexities of the $F_{2i}$, in descending order, and then define a minimum width ordering to be an ordering which gives the least possible sequence with respect to the lexicographic order on sequences of integers. For our purposes, we are only concerned with the maximum complexity of any $F_{2i}$, not how often it occurs. 
This is equivalent to considering the width of the ordering to be the maximum difference between the number of one- and two-handles added at any point in the sequence. We will denote this minimum width over all possible orderings as $c_+(H)$.

The thin position width may just be the same as the genus of the original splitting. For example, a Heegaard splitting is strongly irreducible if any pair of essential discs on opposite sides of the Heegaard surface must intersect, so a strongly irreducible splitting is already in thin position as the order of the handles can not be changed. In \cite{st}, Scharlemann and Thompson show that if a Heegaard splitting is put in thin position, then the manifold is split up in to pieces separated by incompressible surfaces, and each piece contains a strongly irreducible Heegaard splitting. So if the $3$--manifolds $M_i$ have splittings $H_i$ of bounded width $c_+(H_i)$, then each manifold $M_i$ can be split into pieces where each piece has a sweepout of bounded genus. If the Heegaard gradient $\chi_i/d_i$ tends to zero, then the average volume per handle becomes large, so there must be pieces of large diameter, and we can then argue that some of the generalized sweepout surfaces must be parallel.

\subsection{Outline of the proof}

Let $M$ be a closed hyperbolic manifold. Let $M_i$ be a collection of finite covers of $M$, with degree $d_i$, and with Heegaard splittings $H_i$, of Euler characteristic $\chi_i$. The volume of the cover is equal to the volume of $M$ times the degree $d_i$ of the cover, so if the Heegaard gradient $\chi_i/d_i$ tends to zero, then the average volume per handle becomes arbitrarily large. However, we have also assumed that the width of the Scharlemann--Thompson thin position for the Heegaard splittings $H_i$ is bounded, so this means that there is a collection of compression bodies $C_i$ of bounded genus, whose volume tends to infinity.

A sweepout of a $3$--manifold is a one-parameter family of surfaces which fill up the manifold in a degree one manner. A sweepout for a compression body can be constructed by homotoping $\d_+$ down to a spine for the compression body. We can use the hyperbolic metric on the compression body to ``straighten'' the sweepout so that each sweepout surface is negatively curved, except for a single cone point of angle less than $2\pi$. In Section \ref{section:simplicial} we review a precise construction of such sweepouts, called {\sl simplicial sweepouts}, due to Bachman, Cooper and White \cite{bcw}. A technical point is that it is convenient to straighten the sweepout in a complete manifold. We can construct a complete manifold from the compression body by gluing on negatively curved ``flared ends'' to the boundary components, so we end up working in a manifold homeomorphic to the interior of a compression body, with sectional curvature at most $-1$, rather than the original hyperbolic manifold.

One of the main technical tools we will use is a {\sl generalized sweepout}, which is a sweepout in which the genus of the surfaces in the one-parameter family is allowed to change. An ordinary sweepout can be defined as a degree one map from $S \cross I$ to a $3$--manifold, where $S$ is a surface. The domain of the map, $S \cross I$, is a $3$--manifold with height function given by projection onto the $I$ factor. In a generalized sweepout $S \cross I$ is replaced by a $3$--manifold with a Morse function, and the level sets of the Morse function are the sweepout surfaces, which are mapped across into $M$ by a degree one map. We do not require the images of the level sets to be embedded.

The Gauss--Bonnet formula gives an area bound for negatively curved surfaces in terms of their genus, but such surfaces may have arbitrarily large diameter. However, a negatively curved surface has a thick--thin decomposition, so each surface has a ``thick part'', where the injectivity radius is large, and a ``thin part'' where the injectivity radius is small. The thick part has bounded diameter, and the thin part consists of annuli, by the Margulis Lemma. In Section \ref{section:generalized} we give an explicit construction for choosing a continuously varying family of thin parts for the sweepout surfaces, and show how to chop them off to produce a generalised sweepout whose surfaces have bounded diameter.

The sweepout surfaces in the generalized sweepout have bounded diameter, so as the volume of the compression bodies becomes large, there are more and more disjoint sweepout surfaces inside the compression bodies. Although the sweepout surfaces are immersed, it still makes sense to say a sweepout surface $S$ separates a disjoint sweepout surface $S'$ from $\d_+$, if every path from a point on $S'$ to $\d_+$ has $+1$ algebraic intersection number with $S$. We show that as the volume of the compression body becomes large, we can find an arbitrarily large collection of sweepout surfaces which are {\sl nested}, in the sense that for each pair of surfaces, one surface separates the other from $\d_+$. 

We may also assume that the surfaces all have the same genus. As the surfaces are nested, and they are all compressions of $\d_+$, we can show that they are all homotopic. Furthermore, we can find a nested collection of surfaces $S_1, \ldots S_n$ so that the homotopy from $S_n$ to $S_i$ is disjoint from $S_j$ with $j < i$. Gabai \cite{gabai} showed that the singular norm on homology is equal to the Thurston norm, so we can replace the immersed surfaces with embedded surfaces, which we can show are also homotopic. Homotopic surfaces in a compression body need not be isotopic; however if we change the isotopy class of a surface, we must do so by a homotopy that hits a spine of the compression body. As the homotopies between $S_n$ and $S_i$ are disjoint from $S_1$ for $i > 1$, they can't change the isotopy class of the surface, so the surfaces are in fact isotopic, and bound products in the compression body.

So we can find arbitrarily many disjoint parallel embedded surfaces in a compression body of large enough volume. Each compression body lives in a cover of $M$, and each cover is tiled by a choice of fundamental domain for $M$. As the surfaces have bounded diameter, each surface may only hit finitely many fundamental domains of $M$ in the cover $M_i$. So there must be a pair of parallel surfaces which hit the same pattern of fundamental domains, so we can cut the compression body along the matching pair of fundamental domains and then glue them together to form a fibered cover of the original manifold $M$. We only need there to be finitely many ways of gluing the fundamental domains together, so this argument works if we start with a hyperbolic orbifold instead of a manifold.

The remainder of this paper fills in the details of this argument. In Section \ref{section:simplicial} we describe the simplicial sweepout construction due to Bachman, Cooper and White \cite{bcw}, and give a minor generalization to sweepouts of compression bodies. In Section \ref{section:generalized} we describe how to turn a simplicial sweepout of a compression body into a generalized sweepout with sweepout surfaces of bounded diameter. Finally in Section \ref{section:fiber} we complete the proof by showing that as the volume of a compression body tends to infinity we can find arbitrarily many disjoint parallel surfaces of bounded diameter.

\section{Simplicial sweepouts} \label{section:simplicial}

In this section we give a precise definition of a sweepout, and explain how to straighten it to a sweepout with leaves of bounded area.

A {\sl handlebody} is a $3$--manifold with boundary, homeomorphic to a regular neighbourhood of an embedded graph in $\R^3$. A {\sl Heegaard splitting} is a surface which divides a closed $3$--manifold into two handlebodies. A {\sl compression body} is the analogue of a handlebody for $3$--manifolds with boundary, defined as follows. Take $S \cross I$, where $S$ is a closed connected orientable surface, and attach $2$--handles to $S \cross 0$. If any of the resulting boundary components are $2$--spheres, cap them off with $3$--balls. The boundary component corresponding to $S \cross 1$ is called the upper boundary, and will be denoted $\d_+$. The other boundary components are the lower boundary, denoted $\d_-$, which need not be connected. A handlebody is a special case of a compression body in which $\d_-$ is empty.

A {\sl spine} for a handlebody $H$ is an embedded graph $\G$ in $H$ so that the complement $H - \G$ is a product $\d H \cross [0,1)$. Similarly, a spine for a compression body is a $2$--complex $\G$ which is the union of $\d_-$ and a properly embedded graph, so that $C - \G$ is a product $\d_+ \cross [0,1)$. We can construct a sweepout of a handlebody by shrinking the boundary down to the spine, using the product structure. More generally, we shall use the following definition of a sweepout of a compression body.

\begin{definition}[Sweepout]
Let $C$ be a compression body, and let $\G$ be a spine for $C$. Let $S$ be a surface homeomorphic to $\d_+$, and let $\phi\co S \cross I \to C$ be a continuous map. We will write $S_t$ for $S \cross t$. We say the map $\phi$ is a {\sl sweepout} if $\phi$ takes $S_0$ to $\d_- \cup \Gamma$ and $S_1$ to $\d C$, and $\phi$ is degree one, ie, $\phi_*\co  H_3(S \cross I, S \cross \d I) \to H_3(C,\d C)$ is an isomorphism.
\end{definition}

We will often think of the $t$ variable as time. This definition does not require that the surfaces $\phi(S_t)$ be embedded in $C$, and surfaces at different times may intersect each other.

For our purposes, we would like to obtain some control over the geometry of the surfaces $S_t$ in the sweepout. Roughly speaking, we shall do this by triangulating the surfaces, and then ``straightening'' the triangulation in the negatively curved metric on the manifold. The construction we use is that of a {\sl simplicial sweepout}, and is described by Bachman, Cooper and White, in \cite{bcw}. However, we provide a complete definition, as we need a slightly more general construction which allows for manifolds with boundary and change of basepoints.

When we give a surface a triangulation, it need not be a simplicial triangulation, as we wish to allow one-vertex triangulations. Our triangulations will be $\D$--structures as defined by Hatcher \cite{hatcherbook}.

In hyperbolic space the convex hull of three points is a geodesic triangle. In spaces of varying negative curvature this need not be the case, so instead we use coned simplices. A {\sl coned $1$--simplex} $\D^1 = (v_0,v_1)$ is a constant speed geodesic from $v_0$ to $v_1$. We allow degenerate $1$--simplices in which the speed is zero and the image of the $1$--simplex is a point. A {\sl coned $n$--simplex} is a map $\phi\co  \D^n \to M$ so that $\phi | \D^{n-1}$ is a coned $n-1$--simplex, and $\phi | \{ tx + (1-t)v_n | t \in [0,1] \}$ is a constant speed geodesic for all $x \in \D^{n-1}$. Note that the map constructed in this way depends on the order of the vertices in the simplex, and the triangle need not be embedded in $M$.

\begin{definition}[Simplicial surface]
Let $S$ be a triangulated surface, and let $M$ be a closed Riemannian manifold of sectional curvature at most $-1$. Let $\phi\co  S \to M$ be a continuous map so that the map on each triangle $\phi\co \Delta \to M$ is a coned $2$--simplex. Then we say that $\phi\co S \to M$ is a {\sl simplicial surface}.
\end{definition}

A simplicial sweepout is a sweepout in which every sweepout surface is a simplicial surface with a bounded number of triangles, and at most one cone point of positive curvature.

\begin{definition}[Simplicial sweepout]
Let $\phi\co \S \to C$ be a sweepout, such that each surface $S_t$ is mapped to a simplicial surface with at most $4g$ triangles, and with at most one vertex of angle sum less than $2 \pi$. Then we say that $\phi$ is a {\sl simplicial sweepout} of $C$.
\end{definition}

The following theorem follows from the proof of Theorem 2.3 in \cite{bcw}.

\begin{theorem}{\rm\cite[Theorem 2.3]{bcw}}\label{bcw}\qua
Let $M$ be a closed orientable Riemannian manifold of sectional curvature at most $-1$. Let $S_0$ and $S_1$ be simplicial surfaces in $M$, with one vertex triangulations in which the vertices coincide in $M$. Let $\phi\co  S \cross I \to M$ be a homotopy between $S_0$ and $S_1$ which does not move the vertex. Then there is a simplicial sweepout $\Phi'\co  S \cross I \to M$, such that for all $t$, $\Phi(S \cross t)$ consists of at most $4g$ triangles, and has at most one vertex at which the angle sum is less than $2\pi$. Furthermore $\Phi$ is homotopic to $\Phi'$, relative to the boundary of $S \cross I$.
\end{theorem}

In fact we need a minor extension which enables us to change the vertex of the one-vertex triangulation.

\begin{lemma} \label{movebasepoint}
Let $M$ be a closed orientable Riemannian manifold of sectional curvature at most $-1$. Let $S_0$ and $S_1$ simplicial surfaces with one vertex triangulations, which are homotopic by a homotopy $\Phi\co S \cross I \to M$. Then there is a simplicial sweepout $\Phi'\co S \cross I \to M$ which is homotopic to $\Phi$, relative to $S \cross \d I$.
\end{lemma}

\begin{proof}
We reduce to the case of Theorem \ref{bcw} by showing how to homotop $S_0$ so that it shares a common basepoint with $S_1$.

Let $\phi\co  S \cross I \to M$ be the homotopy between the simplicial surfaces $S_0$ and $S_1$. Let $v_0$ and $v_1$ be the basepoints for $S_0$ and $S_1$ respectively. We may assume that the pre-images of the basepoints lie in the same $v \cross I$ fiber in $S \cross I$.  First homotop the homotopy $\phi$ to a new map, which by abuse of notation we will also call $\phi$, in which $\phi$ restricted to $v \cross I$ is a geodesic from $v_0$ to $v_1$. Now homotop $\phi$ to a new map $\phi'$ so that $\phi'(v,t) = \phi(v,2t)$, if $t \leqslant \frac{1}{2}$, and $\phi'(v,t) = \phi(v,1)$ if $t \geqslant \frac{1}{2}$. This means that $\phi'(S \cross \half)$ has the same basepoint as $S_1$, and the homotopy from $S \cross \half$ to $S \cross 1$ is constant on $v \cross t$. Give $S \cross \half$ the same one vertex triangulation as $S \cross 0$, and homotop $\phi'$ so that the edges get mapped to geodesics based at $v_1$. The lower half of $S \cross I$ now has a division into cells which are triangles $\D_i \cross [0,\half]$ in $S \cross [0,\half]$, and $\phi'$ maps the $1$--skeleton of these triangular prisms to geodesic arcs in $M$. So we can homotop $\phi'$ so that $\phi'$ is simplicial on each $\D_i \cross [0,\half]$.

Now $\phi'(S \cross \half)$ and $\phi'(S \cross 1)$ share a common basepoint, and $\phi'$ restricted to $S \cross [\half,1]$ fixes the basepoint, so we can apply Theorem \ref{bcw}.
\end{proof}

\section{Generalised sweepouts} \label{section:generalized}

Simplicial sweepouts are sweepouts in which each surface has bounded area, but we would like to construct sweepouts in which each surface has bounded diameter. In order to do this, we use a more general notion of sweepout, called a generalized sweepout. Recall that a negatively curved surface has a thick--thin decomposition, in which the thick part has bounded diameter, so roughly speaking we wish to cut off the thin part of the sweepout surfaces to obtain bounded diameter surfaces. As the resulting surfaces may become disconnected, it will be convenient to work with the following definition of diameter.

\begin{definition}[$\e$--diameter]
Let $S$ be a set in a metric space $M$. The {\sl $\e$--diameter} of $S$ is the minimum number of balls of radius $\e$ needed to cover $S$.
\end{definition}

The advantage of this definition is that a surface which has two connected components, each of small diameter, will also be considered to have small diameter, even if the two components are far apart.

We now define a generalised sweepout, and give a method for constructing one with bounded $\e$--diameter. Some of the following definitions and constructions are similar to those of Maher and Rubinstein \cite{mr}, but we reproduce them for the convenience of the reader. An important difference here is that the sweepout surfaces are not required to be embedded.

\begin{definition}[Generalised sweepouts]
A {\sl generalised sweepout} of a $3$--man\-if\-old $M$ is a triple $(\S, f, h)$, where

\begin{itemize}

\item $\S$ is an orientable $3$--manifold. 

\item The map $h\co \S \to \R$ is a Morse function, which is constant on each boundary component of $\S$, such that $h^{-1}(t)$ is a collection of surfaces, for all but finitely many $t$.

\item The smooth map $f\co (\S, \d \S) \to (M, \d M)$  is degree one,
  ie, $f_*\co H_3(\S, \d \S)$ $\to H_3(M,\d M)$ is an isomorphism. 

\end{itemize}

We will often write a sweepout as $(\S, \phi)$, where $\phi$ denotes the map $(f \cross h)\co \S \to M \cross \R$.
 We will think of $t \in \R$ as the time coordinate.
\end{definition}

We may think of the sweepout surfaces as a one-parameter family of immersed surfaces in $M$. There may be critical times at which the genus or the number of components of the surface changes. We may change the sweepout in a region $B^3 \cross I \subset M \cross \R$ as follows. Replace $S_t \cap B^3$ with a different continuously varying family of surfaces, so that the new family agrees with the original sweepout on $\d (B^3 \cross I)$. This produces a new generalised sweepout as changing the map in $B^3 \cross I$ doesn't change the degree of the map. We now give a precise definition.

\begin{definition}[Modifying sweepouts]
Let $(\S, \phi)$ be a generalised sweepout. Let $N$ be a $3$--dimensional submanifold of $\S$ whose boundary is disjoint from the critical points of $h$. Let $N'$ be a $3$--manifold whose boundary is homeomorphic to $\d N$, and let $\phi'\co N' \to M \cross \R$ be a map which is the same as the map $\phi$ on $\d N$. Furthermore, assume $\pi_\R \circ \phi'$ is a Morse function on $N'$. Construct a new manifold from $\S$ by cutting out $N$ and replacing it with $N'$, ie, $\S' = (\S - N) \cup N'$. If $\pi_M \circ \phi'\co  \S' \to M$ is still degree one, then $\S'$ is a new generalised sweepout which we say is a {\sl modification} of the original one.
\end{definition}

\begin{figure}[ht!]
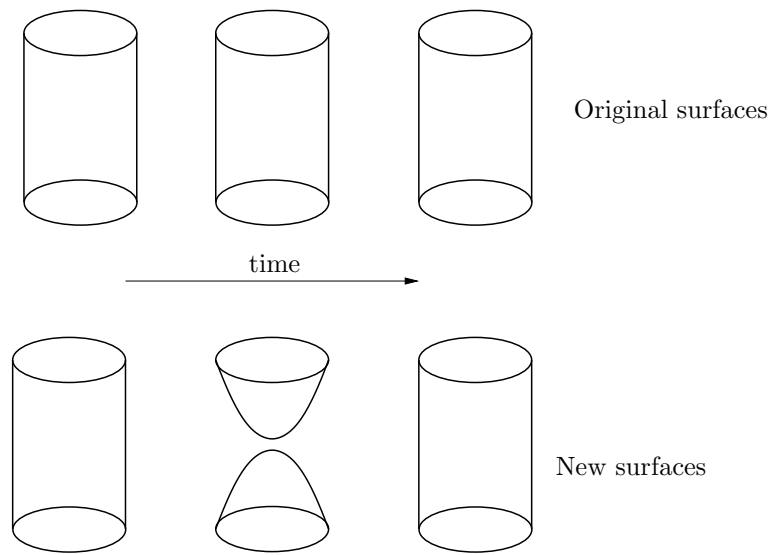
\anchor{picture12}
\begin{center}
\input picture12.pstex_t
\end{center}
\caption{Pinching off a family of annuli}\label{picture12}
\end{figure}

In fact, we will only use one type of modification, which will be to replace a one parameter family of annuli with a family of surfaces which start off as annuli, then get pinched into a pair of discs, and finally get pasted back to the original annulus, as illustrated in \figref{picture12}.

The one parameter family of annuli is a solid torus in $\S$, and this modification corresponds to doing $(0,1)$ surgery on this solid torus. This does not change the degree of the map, as the solid torus represents zero in $H_3(\S,\d \S)$, so the modified map is still a sweepout.

We wish to take the original simplicial sweepout and modify it to produce a generalised sweepout in which each sweepout surface $S_t$ has bounded $\e$--diameter. 

The simplicial sweepout surfaces $S_t$ are composed of triangles of curvature at most $-1$, with at most one vertex of angle sum less than $2\pi$, which we shall call $v_t$. Let $\overline S_t$ be the completion of the universal cover of $S_t - v_t$. This is a simply-connected $2$--complex which is a union of triangles of curvature at most $-1$, with all cone angles at vertices greater than $2\pi$, so it is a complete CAT$(-1)$ geodesic metric space, with the following well known properties, see Bridson and Haefliger \cite{bh}.

\begin{itemize}

\item Geodesics are unique, and hence convex.

\item If $C$ is convex, then $N_r(C)$ is convex.

\item If $C$ is convex, then the closest point projection map onto $C$ is distance decreasing.

\item Isometries are either elliptic, parabolic or hyperbolic. 

\end{itemize}

We remind the reader of the classification of isometries of CAT$(-1)$ spaces. The {\sl translation distance} $D(\phi)$ of an isometry $\phi$ of a metric space $(X,d)$ is $\inf\{ d(x,\phi(x)) | x \in X \}$. The {\sl min set} of $\phi$ is $\{ x \in X | d(x,\phi(x)) = D(\phi)\}$, which may be empty. Elliptic isometries have fixed points, ie, $D(\phi)=0$ and the min set is non-empty, hyperbolic isometries have $D(\phi)>0$ and non-empty min sets, called axes, and parabolic isometries have empty min sets.

For each homotopy class $\a$ in $S_t - v_t$ there is a corresponding covering translation of the universal cover, which gives rise to an isometry of $\overline S_t$. These isometries may be elliptic or hyperbolic, but not parabolic, as the the completion of a preimage of a fundamental domain for $S_t - v_t$ in the universal cover is compact. Let $\g_t$ be the set of points which are moved the least distance by the isometry. These either form a geodesic which is the axis of a hyperbolic isometry, or consist of a single fixed point, if the isometry is elliptic. We can consider $\g_t$ to be the image of the shortest length loop in the homotopy class $\a$. If $\g_t$ consists of a single point, then we consider it to be the constant loop of length zero. If $\g_t$ is a piecewise geodesic associated to $\a$, then it is homotopic to $\a$, if it is disjoint from $v_t$. Otherwise is it is homotopic to $\a$ by an arbitrarily small perturbation which makes it disjoint from $v_t$. By abuse of notation we will call $\g_t$ the geodesic representative of $\a$.

The surfaces $S_t$ are made from negatively curved triangles, which change continuously with $t$, so the geodesic representatives of homotopy classes change continuously with $t$ as well. We now give a detailed proof of this. The basic point is that for surfaces which are close together in time, there is a quasi-isometry between them which is close to an isometry. A geodesic in one surface gets mapped to a quasigeodesic in the other surface, and it is well known \cite{bh} that a quasigeodesic lies in a $K$--neighbourhood of a geodesic for some $K$. It suffices to show that by choosing the times to be sufficiently close together, we may choose $K$ to be as small as we like.

\begin{lemma}Let $\g$ be a simple closed curve in $S - v$. The geodesic representatives $\g_t$ of $\g$ vary continuously with $t$.
\end{lemma}


\begin{proof}
It suffices to consider $S_{t_1}$ and $S_{t_2}$ with $t_1$ and $t_2$ close together. So we may assume there is at most one change of triangulation between $t_1$ and $t_2$. If there is a change of triangulation at time $t_3$, we can break the interval up into two parts $[t_1, t_3]$ and $[t_3, t_2]$, where on each subinterval the surfaces have the same triangulation, with a possibly degenerate metric on some of the triangles. If $\g_t$ varies continuously with $t$ on both closed subintervals, then $\g_t$ varies continuously with $t$ on the whole interval.

Let $\g$ be a simple closed curve in $S_t - v_t$, and let $\g_t$ be the geodesic representative of $\g$ at time $t$.
We want to show that for all $\e > 0$ there is a $\delta > 0$ such that if $| t_1 - t_2  | < \delta$ than $\g_{t_1}$ lies in $N_\e(\g_{t_2})$.

There is a map $\phi\co S_{t_1} \to S_{t_2}$ which is a
$(\l,0)$--quasi-isometry, with $\l \to 1$ as $t_1 \to t_2$. So
$\phi(\g_{t_1})$ is a $(\l,0)$--quasigeodesic in $S_{t_2}$. It is a
standard result of \cite{bh} that therefore $\phi(\g_{t_1})$ lies in a
$N_K(\g_{t_2})$ for some $K$. We wish to show that by choosing $t_1$
close to $t_2$, and hence $\l$ close to $1$, we can make $K$ as small
as we like.

Suppose that $\phi(\g_{t_1})$ is not contained in a $K$--neighbourhood of $\g_{t_2}$. Then there must be a segment of $\phi(\g_{t_1})$ of length at least $K$ lying outside a $K/2$--neighbourhood of $\phi(\g_{t_1})$. Let $x$ and $y$ be the endpoints of this segment, so $d_\phi(x,y) \geqslant K$. Let $d(x,y)$ be the length of the geodesic segment between them in $S_{t_2}$. Let $d_N(x,y)$ be the distance between $x$ and $y$ along the boundary of the $K/2$ neighbourhood. This is illustrated in \figref{picture8}.

\begin{figure}[ht!]
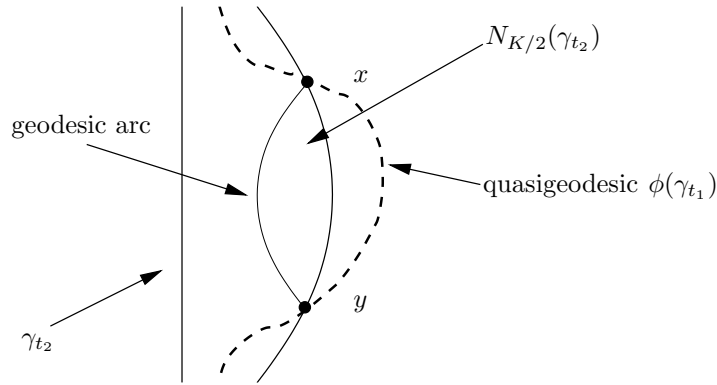
\anchor{picture8}
\begin{center}
\input picture8.pstex_t
\end{center}
\caption{A $K/2$--neighbourhood of a geodesic} \label{picture8}
\end{figure}

A lower bound for the quasigeodesic constant is given by $K/d(x,y)$,
and the largest possible value of $d(x,y)$ occurs when the distance
between $x$ and $y$ along $\d N_{K/2}(\g)$ is $K$. The geodesic
curvature of $\d N_{K/2}(\g)$ is greater than zero if $K>0$, and tends
to zero as $K$ tends to zero. Therefore by choosing $\l$ to be close
to $1$, we can make $K / d(x,y)$ as close to $1$ as we like. This
forces the geodesic curvature to be close to zero, thus making $K$ as
small as we like.
\end{proof}

\begin{definition}
A simple closed curve in $S_t-v_t$ is {\sl short} if the length of its geodesic representative is less than $\e$, the injectivity radius of $M$.
\end{definition}

Let $\Gamma_t$ be the collection of short simple closed curves. There are only finitely many short simple closed curves, and there is always at least one, as the loop around $v_t$ has length zero.

We now define some special neighbourhoods of the short simple closed curves $\gamma_t \in \Gamma_t$ as follows. Consider a connected component $\tilde \g_t$ of the pre-image of $\g_t$ in the completion of the universal cover $\overline S_t$, and choose a transverse orientation for $\tilde \g_t$, so that the distance function from $\tilde \g_t$ is signed. If $\g_t$ has length zero, then assume that distances are positive. For a connected interval $[a,b]$, let $\tilde N_{[a,b]}(\tilde \g_t)$ be the set of points $q$ in $\overline S_t$ with $a \leqslant d(q,\tilde \g_t) \leqslant b$. Let $N_{[a,b]}(\g_t)$ be the projection of $\tilde N_{[a,b]}(\g_t)$ into $S_t$. We will write $N_{[r]}$ if the interval consists of a single point $r$.

\begin{figure}[ht!]
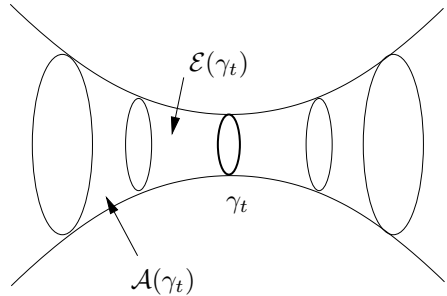
\anchor{picture4}
\begin{center}
\input picture4.pstex_t
\end{center}
\caption{Annular and surgery neighbourhoods}\label{picture4}
\end{figure}

\begin{definition}[Annular and surgery neighbourhoods]
Define an {\sl annular\break neighbourhood} of $\g_t \in \G_t$ as follows. Let $\A(\g_t)$ be the maximal $N_{[a,b]}(\g_t)$ such that each $N_{[r]}(\g_t)$ is embedded of length at most $\e$. The neighbourhood $\A(\g_t)$ is non-empty as it contains $N_{[0]}(\g_t)$. The annular neighbourhoods $\A(\g_t)$ vary continuously with $t$.
See \figref{picture4}.

Given an annular neighbourhood $\A(\g_t) = N_{[a,b]}(\g_t)$ of $\g_t$, define an annular {\sl surgery neighbourhood} $\E(\g_t)$ of $\g_t$ to be the subset of $\A(\g_t)$ of curves $N_{[r]}(\g_t)$ which are at least a distance $\e/2$ from the boundary of $\A(\g_t)$. So $\E(\g_t) = N_{[a+\e/2, b-\e/2]}(\g_t)$, with the convention that this is the empty set if $b-a < \e$. The surgery neighbourhood $\E(\g_t)$ varies continuously with $t$, but need not contain $\g_t$.

\end{definition}

\begin{lemma}
Let $\a_t$ and $\b_t$ be the geodesic representatives of distinct short curves in the sweepout surface $S_t$. Then their surgery neighbourhoods $\E(\a_t)$ and $\E(\b_t)$ are disjoint.
\end{lemma}

\begin{proof}
Suppose $\E(\a_t)$ and $\E(\b_t)$ intersect. Then there is a boundary curve $N_{[r]}(\a_t)$ of $\E(\a_t)$ which intersects $\E(\b_t)$. However the boundary curve $N_{[r]}(\a_t)$ has length at most $\e$, so it must therefore be contained in $\A(\b_t)$. But this means $N_{[r]}(\a_t)$ is either null-homotopic, or homotopic to the other short loop $\b_t$, a contradiction.
\end{proof}

As the surgery neighbourhoods are disjoint, there are at most $2g-1$ such neighbourhoods.

We now show how to create a new generalized sweepout $\hat S_t$ from the original sweepout $S_t$ by cutting off the annular surgery neighbourhoods and replacing them with discs, using the modification described above.

First we describe how to remove a given surgery neighbourhood $\E(\g_t)$. Let $[a,b]$ be a maximal time interval on which $\E(\g_t)$ is non-empty, so $\E(\g_t)$ consists of a single curve $N_{[x]}(\g_a)$ at time $a$, and a single curve $N_{[y]}(\g_b)$ at time $b$. The union of the surgery annuli $\E(\g_{[a,b]})  = \{ \E(\g_t) | t \in [a,b] \}$ is a solid torus in $\S$. We now describe how to produce a new sweepout, identical to the first one outside $\E(\g_{[a,b]})$, and which corresponds to doing $(0,1)$ surgery on $\S$. 

\begin{figure}[ht!]
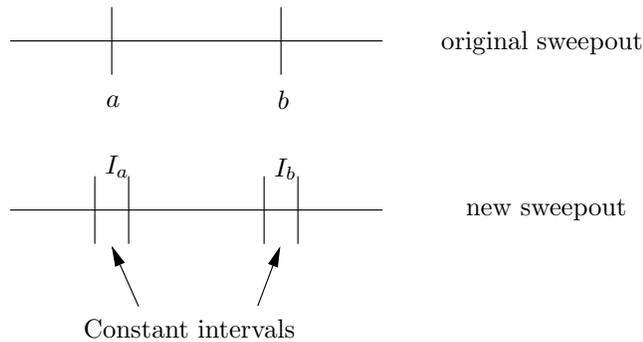
\anchor{picture2}
\begin{center}
\input picture2.pstex_t
\end{center}
\caption{Inserting short constant intervals} \label{picture2}
\end{figure}

Choose a continuously varying family of basepoints for the boundary components of $\E(\g_t)$, which agree at times $a$ and $b$.
Now alter the sweepout by expanding times $a$ and $b$ to short intervals $I_a$ and $I_b$, for which the map $\phi$ is constant, as illustrated in 
\figref{picture2}.

In the interval $I_a$ do a cut move which replaces the loop $\E(\g_a)$ by a pair of discs coned from the basepoint. The coned discs are ruled surfaces, so can be developed into $\H^2$, so their area is bounded above by the area of a disc with perimeter the same length as the curve $\E(\g_a)$, namely $\e$, so each disc has area at most $\e/4$, and is contained in a ball of radius $\e$. In the interval $(a,b)$ replace the annulus $\E(\g_t)$ by the pair of discs formed by coning each boundary component of $\E(\g_t)$ to its basepoint. Finally in the interval $I_b$ paste the two coned discs back to the original surface.

\begin{figure}[ht!]
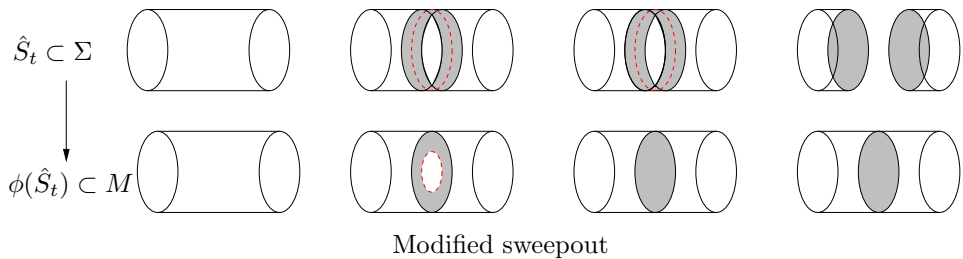
\anchor{picture10}
\begin{center}
\input picture10.pstex_t
\end{center}
\caption{Doing a cut move in $I_a$} \label{picture10}
\end{figure}

We can do this modification for each surgery annulus $\E(\g_{[a,b]})$, and so remove all the surgery annuli from the sweepout, to produce a new sweepout which we shall call $\hat S_t$. As there are at most $2g-1$ short annuli, the new sweepout surface $\hat S_t$ also has area at most $4\pi g+2(2g-1)\e/4$. 

We now recall some useful facts about Voronoi decompositions of surfaces. Let $x_i$ be a collection of points in a metric space $(X,d)$. Let $V_i = \{ x \in X | d(x,x_i) \leqslant d(x,x_j) \text{ for all } j \}$. We say that $V_i$ is the {\sl Voronoi decomposition} of $X$ with respect to $x_i$. If we have chosen a maximal collection of points $x_i$ so that no two are closer together than $\e$, we say we have a {\sl maximal $\e$--spaced Voronoi decomposition}. 

Let $S$ be a surface with a metric, which may be degenerate, so that the induced metric on the universal cover is CAT$(-1)$. It will be useful to know that if we have chosen a maximal $\e$--spaced Voronoi decomposition for $S$, then every Voronoi region which is not a disc has boundary components which are essential in $S$.

\begin{lemma} \label{voronoi_facts}
Let $S$ be a surface with a metric, which may be degenerate, so that the induced metric on the universal cover is CAT$(-1)$, and let $V_i$ be a maximal $\e$--spaced Voronoi decomposition for $S$. Then 

\begin{enumerate}

\item $B_{\e/2}(x_i) \subset V_i \subset B_\e(x_i)$.

\item A Voronoi region which is not a disc has boundary components which are essential in $S$.

\end{enumerate}

\end{lemma}

\begin{proof}
For the first claim, note that any two points $x_i$ and $x_j$ are a distance at least $\e$ apart, so if a point $y$ is at most $\e/2$ from $x_i$, then $y$ must be at least as far from any other $x_j$, so $B_{\e/2}(x_i)$ is contained in $V_i$. If any point in $V_i$ is a distance greater than $\e$ from $x_i$, then it must also be a distance at least $\e$ from all of the other $V_j$, so we could add it to our collection $x_i$ to produce a larger collection of points distance at least $\e$ apart, contradicting maximality.

We now prove the second claim. We have shown that $B_{\e/2}(x_i) \subset V_i \subset B_\e(x_i)$, so the boundary of $V_i$ in $S$ lies in $B_\e(x_i) - B_{\e/2}(x_i)$. If some boundary component of $V_i$ is inessential, then it lifts to a simple closed curve in the universal cover, which lies in $B_\e(\tilde x_i) - B_{\e/2}(\tilde x_i)$, for some pre-image $\tilde x_i$ of $x_i$. This closed curve bounds a disc, and as the distance function $d(x,\tilde x_i)$ is convex, it has maxima only on the boundary, so the entire disc is contained in $B_\e(\tilde x_i)$. However the disc must contain at least one other pre-image of some $x_j$. This can't be a pre-image of $x_i$, as then there would be an essential loop in $S$ intersecting the inessential simple closed once, so it must be a pre-image of some $x_j$ with $j \not = i$. However this implies that there is an $x_i$ and an $x_j$ with $d(x_i,x_j) < \e$, a contradiction.
\end{proof}

Given a surface with a Voronoi decomposition, there is a one-vertex triangulation with edge lengths bounded in terms of the number and size of the Voronoi regions. The surface may have an arbitrary metric.

\begin{lemma} \label{edge_length}
Let $S$ be a closed surface, with a Voronoi decomposition $V_i$ with $N$ regions, such that each $V_i$ is contained in a ball of radius $\e$. Then the surface has a one-vertex triangulation in which each edge has length at most $2\e N$.
\end{lemma}

\begin{proof}
Let $\tilde V_i$ be the Voronoi decomposition of the universal cover $\tilde S$ of $S$, corresponding to the pre-images of the points $x_i$ in $S$. Let $\tilde \G$ be the Delauney graph dual to the Voronoi decomposition of the universal cover. We may choose $\tilde \G$ to be equivariant, and to have edges of length at most $2\e$. Let $\G$ be the projection of $\tilde \G$ into $S$. The complementary regions of $\G$ are all discs, but $\G$ need not be a triangulation, as it may have complementary regions that have more than three sides. 

Complementary regions with more than three sides correspond to points at which more than three distinct $\tilde V_i$ meet. So any pair of pre-images $\tilde x_i$ and $\tilde x_j$ are distance at most $2\e$ apart, so we can add extra edges to triangulate each $n$--gon, to produce a triangulation in which every edge has length at most $2\e$.

\begin{figure}[ht!]
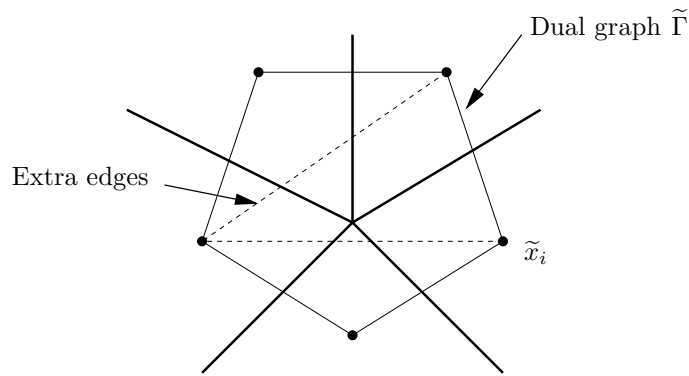
\anchor{picture9}
\begin{center}
\hspace{-2cm}\input picture9.pstex_t
\end{center}
\caption{Add extra edges to make a triangulation} \label{picture9}
\end{figure}

This produces a triangulation of $S$ with the same number of vertices as the number of Voronoi regions, and in which every edge has length at most $2\e$. So we can shrink a maximal tree to a point, and throw away extra edges, to produce a one vertex triangulation in which every edge has length at most $2\e N$.
\end{proof}

\begin{lemma}
There is a constant $N(\e,g)$, which only depends on the genus $g$ of $S_t$ and the injectivity radius $\e$ of the manifold $M$, so that the surface $\hat S_t$ has 

\begin{itemize}

\item $\e$--diameter at most $N$.

\item A one-vertex triangulation in which each edge has length at most $3\e N$.

\end{itemize}

\noindent We may take $N$ to be $12g -10 +4g/\sinh^2(\e/8)$.
\end{lemma}

\begin{proof}
Choose a maximal $\e/2$--spaced Voronoi decomposition for $S_t$, so that $x_1 = v_t$. We now show that there is an upper bound on the number of Voronoi regions. As before, we find bounds on the number of Voronoi regions which are discs, annuli, or have negative Euler characteristic.

If a Voronoi region $V_i$ is a disc, with $i \not = 1$, then by Lemma \ref{voronoi_facts}.1, it contains a ball of radius $\e/4$, so must have area at least $\pi \sinh^2(\e/8)$, so the area of the surface is at most $4\pi g$, there may be at most $4 g /  \sinh^2(\e/8) +1 $ Voronoi regions which are discs.

Voronoi regions with negative Euler characteristic have essential boundary components, so there may be at most $2g-2$ of these. 

If a Voronoi region is an annulus, then it contains a simple closed curve of length at most $\e$. So there are at most $2g-1$ sets of annuli such that the simple closed curves in each set of annuli are parallel. The simple closed curves are at least $\e/2$ apart, so if there are three annuli in a row, the curve corresponding to the middle annulus must be contained in the surgery neighbourhood $\E(\g)$, where $\g$ is the geodesic homotopic to the core curve of the annulus, so there may be at most four annuli for each short curve not contained entirely in the union of the surgery neighbourhoods $\E = \cup \E(\g)$, so there are at most $8(g-1)$ annular Voronoi regions intersecting $S_t - \E$.

Each Voronoi region is contained in a single $\e$--ball in $M$, and we add at most $2g-1$ extra discs to make $\hat S_t$ from $S_t - \E$, each of which is contained in a single $\e$--ball, so the $\e$--diameter of $\hat S_t$ is at most $N = 4g/ \sinh^2(\e/8) +1 +2g-2 +8(g-1) + 2g-1$ as required.

We now show how to construct a triangulation with bounded edge length. Choose a Voronoi decomposition using the points $x_i$ which lie in $S_t - \E$; the above argument shows that there may be at most $N$ such points. Each point in $\hat S_t$ is distance at most $\e/2$ from $S_t - \E$, and hence distance at most $3\e/2$ from an $x_i$ in $S_t - \E$. Therefore the resulting Voronoi decomposition of $\hat S_t$ has each $\hat V_i \subset B_{3\e/2}(x_i)$. The surface $\hat S_t$ now has a triangulation of length $3\e N$ by Lemma \ref{edge_length}.
\end{proof}

\section{Heegaard gradient and virtual fibers} \label{section:fiber}

We briefly summarize of the remainder of the argument to complete the proof of Theorem \ref{main}. If the volume of the covers grows, and $\chi_i/d_i$ tends to zero, but $c_+(M_i)$ remains bounded, then the volumes of the compression bodies in a Scharlemann--Thompson thin position of the manifold must become arbitrarily large. The compression bodies have bounded genus, and hence are swept out by generalised sweepout surfaces of bounded $\e$--diameter. As the sweepout surfaces have bounded $\e$--diameter, they can only intersect finitely many fundamental domains, so eventually there must be many surfaces intersecting the same pattern of fundamental domains in different parts of a cover. Each surface is made from compressing the higher genus boundary of the compression body along discs, and if they are disjoint and nested in the compression body, then they must be homotopic. We use Gabai's result that the singular norm on homology is equal to the Thurston norm to replace the immersed surfaces by embedded surfaces, and then show we can find many parallel embedded surfaces. As the surfaces have bounded diameter, then if there are enough surfaces we can find a pair which intersect matching sets of fundamental domains. We can therefore cut along the matching sets of fundamental domains and glue the ends of the resulting piece of the cover together to form a fibered manifold. 
We now fill in the details of this argument.

If the Heegaard gradient of the covers $\chi_i/d_i$ tends to zero, then the volume of the covers grows faster than the number of handles in the given Heegaard splittings of the covers. The number of compression bodies in the Scharlemann--Thompson thin position for the splitting is bounded by the number of handles, so there must be a sequence of compression bodies in the covers whose volumes become arbitrarily large.

The components of $\d_-$ of a compression body consists of a collection of incompressible surfaces, which may be empty, and the $\d_+$ boundary component consists of a strongly irreducible Heegaard splitting for the union of the two adjacent compression bodies. We need to be able to homotop the boundary to a bounded $\e$--diameter surface close to the original surface, and we will use the fact that the boundary components of a compression body in a thin-position Heegaard untelescoping are closely related to minimal surfaces, which have bounded diameter. 

By Freedman, Hass and Scott \cite{fhs}, the incompressible surfaces are either isotopic to minimal surfaces, or double cover one-sided minimal surfaces in $M$. In the latter case, this would mean that two adjacent compression bodies formed a strongly irreducible splitting of a twisted $I$--bundle over a non-orientable surface. Such a splitting would lift to a Heegaard splitting of the double cover which is an $I$--bundle over an orientable surface, and these splittings are classified by Scharlemann and Thompson \cite{st2}. Given their classification, it is easy to see that the original splitting could not be strongly irreducible. Furthermore, the incompressible surfaces in the untelescoping are all disjoint, so their minimal surface representatives are also all disjoint, as minimal surfaces may only have essential intersections. 

A strongly irreducible Heegaard splitting is either homotopic to a minimal surface, or is homotopic to a minimal surface union some number of arcs, by Pitts and Rubinstein \cite{pr}. We may use the incompressible minimal surfaces as barrier surfaces for the sweepout, so we may assume that the arcs lie in the union of the compression bodies adjacent to the Heegaard splitting surface. As the Heegaard surface is strongly irreducible, the arcs may lie on only one side. If the arcs lie outside the compression body, then we may compress them along discs disjoint from $\d_-$ to produce a new compression body of lower genus, but similar volume. If the arcs lie on the inside of the compression body, we can remove small open neighbourhoods of the arcs so that $\d_+$ still has curvature at most $-1$, while changing the volume of the compression body by an arbitrarily small amount.

We wish to homotop the sweepout to a simplicial sweepout, and we would like to do this in a complete manifold of sectional curvature at most $-1$. We can construct such a manifold by starting with the compression body, and for each boundary component $S$, gluing on a copy of $S \cross [0,\infty)$ with a warped product metric. By the formula for sectional curvature in a warped product (for example see Bishop and O'Neill \cite[page 26]{bo}), given a surface $S$ of Gauss curvature at most $-1$, there is a warped product metric on $S \cross [0,\infty)$ so that the metric is complete, with sectional curvature at most $-1$. This produces a manifold homeomorphic to the interior of the compression body, with a metric of varying negative curvature at most $-1$. We can perturb this metric to make it smooth.

We will show a minimal surface has a triangulation of bounded edge length, where the bound depends on the injectivity radius of the manifold and the genus of the surface. This implies we can straighten the minimal surfaces to simplicial surfaces in a bounded neighbourhood of the original minimal surfaces. If the minimal surface corresponding to the Heegaard surface has extra arcs, then we can triangulate the surface so that any edge in the triangulation runs over the arc at most once, so that the part of the surface corresponding to the arcs collapses down to a geodesic arc, and the rest of the surface has bounded diameter.

At various points we will wish to replace an immersed surface by an embedded surface. An immersed surface $S$ contains embedded surfaces in a regular neighbourhood of $S$ in the same homology class. There is a fundamental result of Gabai \cite{gabai} that the singular norm on homology is the same as the Thurston norm, so we may choose the embedded surface to have genus at most the genus of the immersed surface $S$.

We now show that we can homotop a minimal surfaces to a simplicial surface, which is not too far away in the manifold. This fact follows from the proof of Proposition 6.1 in Lackenby \cite{lac1}, however, for the convenience of the reader, we present a proof, relying only on the following bound from \cite{lac1}.

\begin{lemma}{\rm\cite[proof of Propositon 6.1, claim 1]{lac1}}\label{short homotopy}\qua
Let $F$ be a minimal surface in a Riemannian manifold of curvature at most $-1$, and let $\a$ be a simple closed curve homotopic to a geodesic $\g$. If $\a$ is distance at least $L/2\pi + 1/2$ from $\g$ then the length of $\a$ is at least $L$.
\end{lemma}

We will show that a minimal surface has bounded $\e$--diameter, and a triangulation of bounded edge length.

\begin{lemma} \label{minimal to simplicial}
Let $S$ be a minimal surface in a closed Riemannian manifold $M$ of curvature at most $-1$. Then there is a constant $N$ which depends only on the genus $g$ of $S$ and the injectivity radius $\e$ of $M$, so that $S$ has

\begin{itemize}

\item $\e$--diameter at most $N$.

\item a one-vertex triangulation in which each edge has length at most $2\e N$.

\end{itemize}

\noindent We may take $N$ to be $-\chi(S)(21/4 + 3/4\pi+3/4\e + 2/\sinh^2(\e/4))$. 
\end{lemma}

\begin{proof} Choose a maximal collection of points $x_i$, such that no two points are closer together than $\e$, and let $V_i$ be the Voronoi regions for the $x_i$, ie, $V_i = \{ x \in S | d(x,x_i) \leqslant d(x,x_j) \text{ for all } j \not = i\}$. Each point in a Voronoi region is distance at $\e$ from the basepoint $x_i$ for the Voronoi region, and each Voronoi region contains all the points distance at most $\e/2$ from the basepoint. Each Voronoi region is either a disc, an annulus, or a surface of negative Euler characteristic, and we now consider each of these possibilities in turn.

Each Voronoi region which is not a disc has essential boundary in $S$, by Lemma \ref{voronoi_facts}.2. This means that every Voronoi region of negative Euler characteristic must have boundary components which are essential curves in the surface, and so there may be at most $-3/2\chi(S)$ such pieces.

If a Voronoi region is an annulus, then it contains an essential loop of length at most $\e$, which is homotopic to a geodesic at most a distance $\e/2\pi + 1/2$ away, by Lemma \ref{short homotopy}. This means there may be at most $1/2\e + 5/2 + 1/2\pi$ annuli parallel to any given curve, so there may be at most $-3/2\chi(S)(1/2\e + 5/2 + 1/2\pi)$ annuli in total.

If a Voronoi region is a disc, then it contains a ball of radius $\e/2$ in the surface. By the monotonicity formula for minimal surfaces in $3$--manifolds of negative curvature at most $-1$, for example, see Choe \cite{choe}, the area of this disc must be at least $\pi \sinh^2(\e/4)$, so there are at most $-2\chi(S)/\sinh^2(\e/4)$ discs.

So $N = -\chi(S)(21/4 + 3/4\pi+3/4\e + 2/\sinh^2(\e/4))$ is an upper bound for the total number of Voronoi regions, and hence an upper bound for the $\e$--diameter of the surface. Then by Lemma \ref{edge_length}, the surface has a triangulation in which each edge has length at most $2\e N$.
\end{proof}

Given simplicial surfaces homotopic to the components of $\d_-$, we can join them together with arcs to form a spine for the compression body. We next show that we can find a simplicial surface homotopic to this spine, by a homotopy which does not sweep out very much volume.

\begin{lemma}
Let $S_1, \ldots, S_n$ be a collection of simplicial surfaces, with basepoints $v_i$, and of total genus $g$ in a complete Riemannian $3$--manifold of curvature at most $-1$. Connect the basepoint $v_1$ to each of the other basepoints by a geodesic arc to form a $2$--complex $\S$. Then there is a simplicial surface of genus $g$ which is homotopic to $\S$, by a homotopy which sweeps out a volume of at most $3(2g+2)$ times the maximal volume of an ideal hyperbolic tetrahedron.
\end{lemma}

\begin{proof}
Use the construction from Lemma \ref{movebasepoint} to homotop the basepoint of each $S_i$ to $v_1$, for $i \not = 1$. Each simplicial triangle sweeps out a triangular prism, which may not be embedded. However, the pullback metric on the triangular prism in the domain of the sweepout has curvature at most $-1$, and so is CAT($-1$). Therefore there is a comparison geodesic triangular prism in hyperbolic space, which has the property that distances between points in the hyperbolic metric are at least as big as distances between points in the pullback metric. As a triangular prism can be triangulated with three tetrahedra, the volume of the prism in the pullback metric is at most three times the maximal volume of an ideal hyperbolic tetrahedron. So each triangle sweeps out a volume at most $3$ times $\D$, the maximal volume of an ideal tetrahedron, so the total homotopy sweeps out a volume of at most $3(2g+2)\D$, as there are at most $2g+2$ triangles. Now the surfaces all share a common basepoint, so we can map a simplicial surface of genus $g$ with a one vertex triangulation onto them, so that the vertex goes to the basepoint.
\end{proof}

These two lemmas show that if we have a compression body in a complete Riemannian manifold of curvature at most $-1$, we can construct a sweepout of the compression body, which is simplicial, except possibly in a bounded neighbourhood of the boundary, where the bound depends only on the genus of the boundary.

\begin{definition}
We say a collection of surfaces in a compression body is {\sl nested} if for any pair of surfaces, one of them separates the other from the higher genus boundary of the compression body. 
\end{definition}

As the volume of a compression body $C$ becomes large, the diameter of $C$ must also become large. We use this to show there is a cover which contains a compression body with many disjoint nested sweepout surfaces.

\begin{lemma}
If the volume of a compression body $C$ in one of the covers $M_i$ is at least vol($B_{Kn}(x)$), then there are at least $n$ disjoint nested surfaces, each one of which is formed by compressing the higher genus boundary of the compression body. The constant $K$ depends on $\e$ and the genus of the compression body.
\end{lemma}

\begin{proof}
There is a generalised sweepout $S_t$ of the compression body such that every sweepout surface has $\e$--diameter at most $K(g,\e)$, where $g$ is the genus of the handlebody, and $\e$ is the injectivity radius of $M$. Let $A$ be an upper bound for the volume of a $K\e$--neighbourhood of a sweepout surface $S_t$. Let $A=K \vol (B_{K\e + \e})$.

If the diameter of $C$ is $d$, then there is a point $x$ of distance at least $d/2 - K$ from $\d_+$. Let $\g$ be the geodesic arc from $x$ to $\d_+$. Let $\g$ have length $L$, which is at least $d/2-K$, and let $\g(l)$ be an arc-length parameterisation of $\g$, so that $\g(0)$ is in $\d_+$.

The sweepout surfaces $S_t$ sweep out the compression body $C$. More precisely, we say that a point $x \in C - S_t$ is separated from $\d_+$ by $S_t$ if any path from $x$ to $\d_+$ has algebraic intersection number $+1$ with $S_t$. Let $C_t$ be the subset of $C$ swept out by $S_t$, ie, $C_t$ is the closure of $ \{x \in C | x \text{ is separated from } \d_+ \text{ by } S_t \}$. Then $C_t$ starts out at as the empty set, and ends up as all of $C$. Furthermore, $C_t$ varies continuously with $t$, as $S_t$ varies continuously with $t$. In particular, any subset of $C$, such as $\g$, is also swept out by the sweepout $S_t$, ie, $\g \cap C_t$ varies continuously with $S_t$, and $\g \cap C_0$ is the empty set, and $\g \cap C_1$ is all of $\g$.

There is a sweepout surface $S_t$ which divides $\g$ into two parts of equal length, ie, the length of $C_t \cap \g$ is half the length of $\g$. The sweepout surface $S_t$ may have at most $g$ components. Each component of $S_t$ has diameter at most $K$, so if $\g$ hits some component of $S_t$, and then travels a distance greater than $K$ past that component of $S_t$, then it can never return and hit the same component of $S_t$ again. As $S_t$ splits $\g$ in two, there most be some component of $S_t$ which intersects $\{ \g(l) | l \in ([L/4,3L/4]) \}$ which separates $\g$ into two pieces each of length at least $L/4$. Choose this component of $S_t$ to be $T_1$. Note that $T_1$ separates $x$ from $\d_+$, is disjoint from both of them, and is distance at least $d/8-K$ from $x$. Let $l_1$ be the smallest value of $l$ such that $\g(l) \in T_1$.

So we can now apply the same argument again, this time using $\{ \g(l) | l \in [0,l_1] \}$ instead of $\g$. We can continue in this way finding a sequence of $n$ surfaces $T_i$, until $L/2^{n+1}$ is comparable to $K$.
\end{proof}

Therefore, as the volume of a compression body becomes large, we can find many disjoint nested surfaces inside it, formed from compressing the higher genus boundary, so in particular we can find many nested connected surfaces of the same genus.

We next show that if two connected sweepout surfaces are nested in the compression body, then the compressing discs for the outer one may be chosen to be a subset of the compressing discs for the inner one. In particular, if the surfaces have the same genus, then the compressing discs are the same, so they are homotopic.

\begin{lemma}
Let $S_1$ and $S_2$ be two disjoint immersed sweepout surfaces each of which is homotopic to a surface obtained from compressing $\d_+$ along collections of discs $\D_1$ and $\D_2$ respectively. Furthermore suppose that $S_2$ separates $S_1$ from $\d_+$. Then we may choose compressing discs $\D'_1$ and $\D'_2$ yielding $S_1$ and $S_2$, so that the discs $\D'_2$ are a subset of the discs $\D'_1$.
\end{lemma}

\begin{proof}
The sweepout gives a homotopy from $S_1$ to $\d_+ \cup \D_1$. We may homotop the compressing discs $\D_1$ so they are disjoint and embedded in the compression body $C$. Similarly we can homotop the compressing discs $\D_2$ so they are embedded and disjoint. Note that two different sets of compressing discs $\D_1$ and $\D'_1$ can give rise to isotopic surfaces.

We say a properly embedded essential disc $D$ in $C$ is a compressing disc for $S_i$ if there is a collection of disjoint compressing discs $\D'_i$ containing $D$, so that $\d_+$ cut along $\D'_i$ is isotopic to $\d_+$ cut along $\D_i$.

We first prove the lemma for non-separating discs. Let $D$ be a non-separating essential disc in $C$, which is \emph{not} a compressing disc for $S_1$. Then there is a curve $\a$ in $\d_+ - \D_1$, which has algebraic intersection number one with $D$, and this is preserved under homotopy, so $S_1$ must intersect $D$. As $S_2$ separates $S_1$ from $\d_+$, $S_2$ must also intersect $D$, and in fact the intersection of $S_2$ with $D$ must separate the intersections of $S_1$ in $D$ from $\d D$. Furthermore, there must be an essential intersection of $D$ with $S_2$, as otherwise discs in $S_2$ would have to intersect the image of $\a$ in $S_1$, a contradiction. This means that $D$ may not be a compressing disc for $S_2$ either.

If $D$ is a separating compressing disc for $C$, which is \emph{not} a compressing disc for $S_1$, then it still has homotopically essential intersection with $S_1$. Explicitly, $\D$ separates $\d_+$ into two components, neither of which is planar. Choose curves $\a$ and $\beta$ in each component which are essential in $C$, and connect them with an arc $\g$ which intersects $D$ once. The resulting trivalent graph, or eyeglass, cannot be homotoped to be disjoint from $D$. This can be seen by considering a connected pre-image $\theta$ in the universal cover, which intersects a particular lift $\tilde D$ of $D$. If the pre-image of the graph could be homotoped off $\tilde D$, then this homotopy would lift to a homotopy of the pre-image in the universal cover, but $\theta$ has unbounded components on both sides of $\tilde D$, and the homotopy can only move each point a finite distance, a contradiction. So $S_1$ intersects $D$, as must $S_2$, as $S_2$ separates $S_1$ from $\d_+$. There must be essential intersections of $S_2$ with $D$, as if all intersections were inessential then in the universal cover an unbounded part of $\theta$ would be separated from the pre-image of $\d_+$ by a bounded part of $S_2$, a contradiction. So $D$ may not be a compressing disc for $S_2$.

So we have shown that a compressing disc for $S_2$ is also a compressing disc for $S_1$, so we may choose a set of compressing discs for $S_2$, and then extend this to a set of compressing discs for $S_1$.
\end{proof}

Ultimately, we will need to obtain embedded surfaces, rather than immersed surfaces. As the singular norm on homology is equal to the Thurston norm on homology, we can replace an immersed surface $S_i$ with an embedded surface contained in a regular neighbourhood of $S_i$, in the same homology class, and with genus at most the genus of $S_i$. We now show that if an embedded surface separates one of the surfaces $S_i$ from $\d_+$, then the genus of the embedded surface must be at least as large as the genus of the surface $S_i$.

\begin{lemma} \label{embedded_genus}
Let $S$ be a connected immersed surface in a compression body, which is homotopic to $\d_+$ compressed along some collection of discs. Let $T$ be a least genus connected embedded surface which separates $S$ from $\d_+$. Then $T$ is incompressible in $C - S$, and the genus of $T$ is at least as big as the genus of $S$.
\end{lemma}

\begin{proof}
As $T$ is an embedded surface in a compression body, it is separating, and it can be compressed to a surface parallel to components of $\d_-$ by a sequence of compressions along embedded discs. If any compressing disc is disjoint from $S$ then we can compress along the disc to reduce the genus of $T$, so $T$ is incompressible in $C - S$.

The surface $S$ is homotopic to $\d_+$ compressed along some collection of discs $D$ say, so $C - D$ is a compression body $C'$, and we can choose a spine $\G$ for $C'$ which can be homotoped into $S$. The map on first homology induced by inclusion $H_1(\G) \to H_1(C)$ is injective. As $T$ bounds discs on one side only, $T$ bounds a compression body $C''$ in $C$. Consider the maps induced by inclusion $H_1(\G) \to H_1(C'') \to H_1(C)$. This composition must be injective, so the rank of $H_1(C'')$ is at least as big as the rank of $H_1(\G)$, so by Poincar\'e--Lefschetz duality, the rank of $H_1(T)$ is at least as big as the rank of $H_1(S)$, which implies the genus of $T$ is at least as big as the genus of $S$, as required.
\end{proof}

We have shown that if the volume of a compression body is sufficiently large, there must be many disjoint nested homotopic surfaces inside it. However, we require the further property that the homotopy from $S_n$ to $S_i$ is disjoint from $S_j$ for $j < i$. We now show that we can find a collection of surfaces with this additional property.

\begin{lemma}
Let $S_1, \ldots S_n$ be a collection of disjoint, nested, homotopic surfaces in a compression body $C$. Then there is a collection of surfaces\break $S'_1, \ldots, S'_{n-1}, S_n$ which are disjoint, nested and homotopic, and furthermore the homotopy from $S_n$ to $S'_i$ is disjoint from $S'_j$ for $j < i$. 
\end{lemma}

\begin{proof}
Each surface $S_i$ has a triangulation of bounded length, so the surface is homotopic to a simplicial surface of bounded $\e$--diameter, a bounded distance away in $M$. Therefore there is a homotopy which is a simplicial sweepout turned into a bounded diameter sweepout $S'_t$ from $S_n$ to $S_1$. Let $S'_i$ be the first $S'_t$ which hits $S_i$, so the sweepout from $S_n$ to $S'_i$ is disjoint from $S'_j$ for $j < i$ by construction. Suppose the genus of some $S'_i$ is less than the genus of $S_i$. As the Thurston norm is equal to the singular norm, there is therefore an embedded surface $T'_i$ of genus lower than $S_i$ separating $S_1$ from $S_n$, contradicting Lemma \ref{embedded_genus}. So in fact no compressions take place, and the $S'_i$ are all homotopic, and the homotopy from $S_n$ to $S'_i$ is disjoint from $S'_j$ for $j < i$. The $S'_i$ are disjoint as they have bounded diameter, and are far apart. Finally, we show that $S'_1, \ldots S'_{n-1}, S_n$ are nested. Let $C_n$ be the subset of $C$ separated from $\d_+$ by $S_n$. First note that $S'_{n-1}$ intersects $S_{n-1}$, which is far from $S_n$, and hence far from $\d C_n$, so $S'_{n-1} \subset C_n$, so $S'_{n-1}$ is nested with respect to $S_n$. Similarly, $S'_{n-k}$ is contained within $C_n$, so for any point $x$ in $C_{n-k}$ there is a path $\g$ from $x$ to $\d_+$ with algebraic intersection number $1$ with $S_n$. Now $S'_{n-k+1}$ is homotopic to $S_n$ by a homotopy which is disjoint from $C'_{n-k}$, so this does not change the algebraic intersection number, so $S'_{n-k} \subset C'_{n-k+1}$, as required.
\end{proof}

If we choose $S$ to be $S_1$, and $T$ to be an embedded surface corresponding to $S_i$, with $i > 1$, then Lemma \ref{embedded_genus} shows that the genus of the embedded surface is the same as genus of the immersed surface $S_i$, and furthermore, that the embedded surface is incompressible in the complement of $S_1$.

Any two embedded surfaces $T_i$ and $T_j$, with $1 < i < j < n$ bound a submanifold in the compression body $C$, and we now show that any homotopy between $S_1$ and $S_n$ sweeps out this submanifold in a degree one manner.

\begin{lemma}
Let $S_1, \ldots, S_n$ be a collection of disjoint connected nested immersed homotopic surfaces in a compression body $C$, and let $T_i \subset N(S_i)$ be embedded surfaces in the same homology class. Let $X$ be the part of $C$ between $T_{n-1}$ and $T_2$. Then any homotopy $\phi\co S \cross I \to C$ from $S_1$ to $S_n$ is degree one onto $X$, ie, $\phi_*\co H_3(S \cross I, \d) \to H_3(X,\d X)$ is an isomorphism. 
\end{lemma}

\begin{proof}
The map $\phi\co S \cross I \to C$ gives a map of pairs $\phi\co (S \cross I,\d) \to (C,C-X)$, and by excision, $H_3(C, C - X_2) \cong H_3(X, \d X)$. Continuous maps between pairs induce natural homomorphisms between the homology long exact sequences of pairs, so the following diagram commutes.
\begin{equation*}
\begin{CD}
H_{3}( S \cross I,\d)  @>\d>> H_{2}(S \cross 0 \cup S \cross 1) \\
@V \phi_* VV  @V \phi_* VV \\
H_{3}(C, C - X)  @> \d >> H_2(C-X)  \\
\end{CD}
\end{equation*}
The group $H_3(S \cross I, \d)$ is generated by $[S \cross I]$, which is mapped to $[S \cross 1] - [S \cross 0]$ by the boundary map, and $\phi_*$ maps this to $[S_n] - [S_1]$, which is homologous to $[T_{n-1}] - [T_2]$. If there is no component of $\d_-$ contained in the component of $C-X$ containing $T_2$, then $[T_2]$ will be zero. However, $[T_{n-2}]$ is not zero, as it lies in the same component of $C - X$ as $\d_+$.

The group $H_{3}(C, C - X) \cong H_3(X, \d X)$ is generated by $[X]$, which gets mapped to $[T_{n-1}] - [T_2]$ by the boundary map, so $[S \cross I]$ must get mapped to $[X]$, as required.
\end{proof}

We next show that the embedded surfaces are in fact homotopic to the original immersed surfaces.

\begin{lemma}
Let $S_1$ and $S_n$ be disjoint immersed surfaces in a compression body $C$, which are homotopic, by a homotopy $\phi\co S \cross I \to C$. Let $T$ be an embedded surface which separates $S_1$ from $S_n$, which is incompressible in $C - (S_1 \cup S_2)$, and for which the map $\phi_*\co H_3(S \cross I, \d) \to H_3(N(T),\d)$ is degree one, where $N(T)$ is a regular neighbourhood of $T$. Then $T$ is homotopic to $S_t$.
\end{lemma}

\begin{proof}
Let $\phi\co S \cross I \to M$ be a smooth homotopy from $S_1$ to $S_n$. The embedded surface $T$ is contained in the image of $\phi$, and $T$ separates $S_1$ from $S_n$, so $\phi^{-1}(T)$ is an embedded surface in $S \cross I$ which separates $S \cross 0$ from $S \cross 1$. Furthermore the map $\phi$ is degree one from $\phi^{-1}(T)$ to $T$. 

Suppose the surface $\phi^{-1}(T)$ is compressible in $S \cross I$. By the loop theorem, we may assume the compressing discs are embedded. Let $D$ be a maximal collection of disjoint embedded compressing discs for $\phi^{-1}(T)$. Construct a map of a surface $\psi \co  F \to \phi^{-1}(T) \cup D$ as follows. Take the closure of $\phi^{-1}(T) - D$ with the induced path metric to produce a surface $F'$ with $2|D|$ boundary components. Now glue in a disc along each boundary component, and let $F$ be the resulting surface. We now define a map from $F$ to $\phi^{-1}(T) \cup D$. Send each point in $F'$ to the corresponding point in $\phi^{-1}(T)$. Send each disc in $F - F'$ to the corresponding compressing disc in $S \cross I$. This is illustrated in \figref{picture6}.

\begin{figure}[ht!]
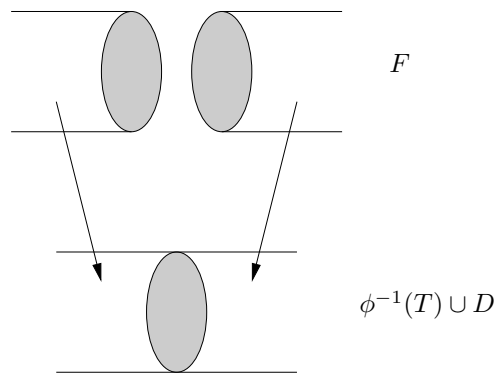
\anchor{picture6}
\begin{center}
\input picture13.pstex_t
\end{center}
\caption{A map from $F$ to $\phi^{-1}(T) \cup D$} \label{picture6}
\end{figure}

The surface $T$ is incompressible in the complement of $S_1$, so the
images of the discs $D$ in $M$ are all inessential, and so they are
homotopic into $T$, as $M$ is irreducible. So the map $\phi \circ
\psi\co F \to M$ is homotopic onto the image of $T$, and is still
degree one, as if it were degree zero, then there would be a path in
$C$ from $S_1$ to $S_n$ with zero algebraic intersection number with
$\phi \circ \psi (F)$.
This path would pull back to a path from $S \cross 0$ to $S \cross 1$
with zero intersection number with $\psi (F)$, but $\psi (F)$
separates $S \cross 0$ from $S \cross 1$, a contradiction. This means
that $F$ is in fact homotopic to $T$ in $C$. Furthermore, $\psi(F)$ is
an incompressible surface in $S \cross I$, so is homotopic to $S
\cross I$, hence $\phi \circ \psi(F)$ is homotopic to $S_t$ in $C$, as
required.
\end{proof}

We now show that the embedded surfaces $T_i$ are not just homotopic, but actually parallel in the compression body $C$.

\begin{lemma}
Let $T_{n-1}$ and $T_3$ be embedded incompressible surfaces contained in the image of a homotopy $\phi$ from $S_n$ to $S'_2$, such that the image of $\phi$ is disjoint from $S'_1$. Then $T_{n-1}$ and $T_3$ are parallel.
\end{lemma}

\begin{proof}
Let $Y$ be the region bounded by $T_{n-1}$ and $T_3$. Consider $\phi^{-1}(Y)$ in $S \cross I$, \figref{picture7}.

\begin{figure}[ht!]\anchor{picture7}
\begin{center}
\input picture7.pstex_t
\end{center}
\caption{$\phi^{-1}(Y)$} \label{picture7}
\end{figure}

Suppose some component of $\phi^{-1}(Y)$ has compressible boundary. Then we can change the homotopy $\phi$ in a neighbourhood of the compressing disc to compress the boundary. If any component of $\phi^{-1}(Y)$ is a $3$--ball, we can change the homotopy in a neighbourhood of the $3$--ball to remove it. The remaining components of $\phi^{-1}(Y)$ must all have boundary components isotopic to $S \cross t$. We can change the product structure on $S \cross I$ so that these boundary components are horizontal, and then there is a subset of $I$ which gives a homotopy from $T_{n-1}$ to $T_3$ entirely contained in $Y$. So now $Y$ is a $3$--manifold with two incompressible boundary components, which are homotopic inside $Y$, so $Y$ has the same fundamental group as a surface, so $Y$ is a product $S \cross I$, by Waldhausen \cite{waldhausen}.
\end{proof}

\begin{lemma}
If there are enough disjoint parallel surfaces then the manifold is virtually fibered.
\end{lemma}

\begin{proof}
A choice of fundamental domain for the original orbifold gives a tiling of any cover. Each parallel surface has bounded diameter, so is contained in finitely many fundamental domains, and we can choose parallel surfaces so that the fundamental domains they hit are disjoint. There are only finitely many ways of gluing finitely many fundamental domains together, so if there are enough parallel surfaces, there must be at least three which hit the same pattern of fundamental domains, so we cut the manifold along a pair of fundamental domains with compatible orientation, and glue back together to get a fibered cover of the original orbifold $X$. 
\end{proof}

This completes the proof of Theorem \ref{main}.


\end{document}

%% file: gtoutput.tex
\def\ifplaintex{\expandafter\ifx\csname documentclass\endcsname\relax}

\ifplaintex 
\else
\fi

\expandafter\ifx\csname beginpicture\endcsname\relax
\expandafter\ifx\csname documentclass\endcsname\relax
\input pictex \else
\input prepictex \input pictex \input postpictex \fi\fi

\def\gt{{\mathsurround=0pt\it $\cal G\mskip-2mu$eometry \&\ 
$\cal T\!\!$opology}}        

\def\gtp{{\mathsurround=0pt\it $\cal G\mskip-2mu$eometry \&\ 
$\cal T\!\!$opology $\cal P\!$ublications}}  

\def\lognumber#1{\def\thelognumber{#1}}
\def\volumenumber#1{\def\thevolumenumber{#1}}
\def\papernumber#1{\def\thepapernumber{#1}}
\def\volumeyear#1{\def\thevolumeyear{#1}}

\def\pagenumbers#1#2{\def\startpage{#1}\def\finishpage{#2}}
\def\published#1{\def\publishdate{#1}}
\def\proposed#1{\def\theproposer{#1}}
\def\seconded#1{\def\theseconders{#1}}
\def\received#1{\def\receiveddate{#1}}

\def\accepted#1{\def\accepteddate{#1}}

\long\def\asciiabstract#1{\long\def\theasciiabstract{#1}}
\def\asciikeywords#1{\def\theasciikeywords{#1}}

\let\\\par\let\thelognumber\relax
\let\thevolumenumber\relax\let\thepapernumber\relax
\let\thevolumeyear\relax\let\thesamplenumber\relax\let\startpage\relax
\let\finishpage\relax\let\publishdate\relax\let\receiveddate\relax
\let\reviseddate\relax\let\accepteddate\relax\let\theasciititle\relax
\let\theasciiauthors\relax
\let\theasciiabstract\relax\let\theasciikeywords\relax
\let\theasciiemail\relax\let\theshortauthors\relax\let\theshorttitle\relax

\long\def\maketitlep{   

\count0=\startpage

\gt\hfill      
\beginpicture
\setcoordinatesystem units <0.33truein, 0.33truein> point at 2.2 0.9
\setplotsymbol ({$\cal G$})
\plotsymbolspacing=9truept
\circulararc 315 degrees from 0 1 center at 0 0
\setplotsymbol ({$\cal T$})
\circulararc 315 degrees from 1 -1 center at 1 0
\endpicture
\break
{\small\ifx\thesamplenumber\relax 
Volume \else Sample
\fi\thevolumenumber\ (\thevolumeyear)
\startpage--\finishpage\nl
Published: \publishdate}
\vglue 0.5truein plus 0.4fil minus 0.1truein

{\parskip=0pt\leftskip 0pt plus 1fil\def\\{\par\smallskip}{\ifplaintex\large
\else\Large\fi\bf\thetitle}\par\medskip}   

\vglue 0pt plus 0.1fil 

{\parskip=0pt\leftskip 0pt plus 1fil\def\\{\par}{\sc\theauthors}
\par\medskip}

\vglue 0pt plus 0.1fil 

{\small\parskip=0pt\let\newline\\
{\leftskip 0pt plus 1fil\def\\{\par}{\sl\theaddress}\par}
\expandafter\ifx\theemail\relax    
\relax\else\vglue 5pt plus 0.02fil minus 2pt\def\\{\stdspace{\rm 
and}\stdspace} 
\cl{Email:\stdspace\tt\theemail}\fi
\ifx\theurl\relax                  
\relax\else\vglue 5pt plus 0.02fil minus 2pt\def\\{\stdspace{\rm 
and}\stdspace}
\cl{URL:\stdspace\tt\theurl}\fi\par}

\vglue 7pt plus 0.3fil minus 3pt

{\bf Abstract}
\vglue 5pt plus 0.1fil minus 2pt

\theabstract

\vglue 7pt plus 0.3fil minus 3pt

{\bf AMS Classification numbers}\quad Primary:\quad \theprimaryclass

Secondary:\quad \thesecondaryclass

\vglue 5pt plus 0.3fil minus 2pt

{\bf Keywords:}\quad \thekeywords

\vglue 10pt plus 0.5fil minus 5pt

{\small  Proposed: \theproposer\hfill Received: \receiveddate\nl
Seconded: \theseconders\hfill 
\ifx\reviseddate\relax                         
Accepted: \accepteddate                        
\else
Revised: \reviseddate                          
\fi}
\eject
}       

\let\maketitlepage\maketitlep
\let\maketitle\maketitlepage

\font\phead=cmsl9 scaled 950
\font=cmsl9 scaled 1050
\font\pnum=cmbx10 scaled 913
\font\lnum=cmbx10 
\font\pfoot=cmsl9 scaled 950
\font=cmsl9 scaled 1050
\ifplaintex
\headline{\vbox to 0pt{\vskip -4.5mm\line{\small\phead\ifnum
\count0=\startpage ISSN 1364-0380 (on line)
1465-3060 (printed) \hfill {\pnum\folio}\else\ifodd\count0\def\\{ }%
\ifx\theshorttitle\relax\thetitle\else\theshorttitle\fi\hfill{\pnum\folio}
\else\def\\{ and }{\pnum\folio}\hfill\ifx\theshortauthors\relax\theauthors
\else\theshortauthors\fi\fi\fi}\vss}}
\footline{\vbox to 0pt{\vglue 0mm\line{\small\pfoot\ifnum\count0=\startpage
\copyright\ \gtp\hfill\else
\gt, Volume \thevolumenumber\ (\thevolumeyear)\hfill\fi}\vss
}}
\else
\makeatletter
\def\@oddhead{{\small\lhead\ifnum\count0=\startpage ISSN 1364-0380 (on line)
1465-3060 (printed) \hfill {\lnum\number\count0}\else\ifodd\count0
\def\\{ }\ifx\theshorttitle\relax \thetitle \else\theshorttitle\fi\hfill
{\lnum\number\count0}\else\def\\{ and }{\lnum\number\count0}
\hfill\ifx\theshortauthors\relax 
\theauthors\else\theshortauthors\fi\fi\fi}}\def\@evenhead{\@oddhead}
\def\@oddfoot{\small\lfoot\ifnum\count0=\startpage\copyright\ \gtp\hfill\else
\gt, Volume \thevolumenumber\ (\thevolumeyear)\hfill\fi}
\def\@evenfoot{\@oddfoot}
\makeatother
\fi

\newwrite\gtoutfile
\long\gdef\makeheadfile{  
{\def\\{, }\def\s{ }
\immediate\openout\gtoutfile head.xxx
\immediate\write\gtoutfile{Proxy-for: \ifx\theasciiauthors\relax
\theauthors\else\theasciiauthors\fi\s<\ifx\theasciiemail\relax\theemail\else\theasciiemail\fi>}
\immediate\write\gtoutfile{\noexpand\\}
\immediate\write\gtoutfile{Authors: \ifx\theasciiauthors\relax
\theauthors\else\theasciiauthors\fi}
{\def\\{ }\immediate\write\gtoutfile{Title: \ifx\theasciititle\relax
\thetitle\else\theasciititle\fi}}
\immediate\write\gtoutfile{Subj-class: GT or SG or MG etc}
\immediate\write\gtoutfile{MSC-class: \theprimaryclass\ifx\thesecondaryclass\relax\else, \thesecondaryclass\fi}
\immediate\write\gtoutfile{Journal-ref: Geom. Topol. \thevolumenumber
(\thevolumeyear) \startpage-\finishpage}
\immediate\write\gtoutfile{Comments: Published by Geometry and Topology at}
\immediate\write\gtoutfile{\s\s http://www.maths.warwick.ac.uk/gt/GTVol\thevolumenumber/paper\thepapernumber.abs.html}
\immediate\write\gtoutfile{\noexpand\\}
\immediate\write\gtoutfile{}
\ifx\theasciiabstract\relax
\immediate\write\gtoutfile{\theabstract}\else
\immediate\write\gtoutfile{\theasciiabstract}\fi
\immediate\write\gtoutfile{}
\immediate\write\gtoutfile{\noexpand\\}
\immediate\write\gtoutfile{}
\immediate\closeout\gtoutfile}}  

\def\maketitlepage{\maketitlep\makeheadfile}
\let\maketitle\maketitlepage

%% file: picture12.pstex_t
\begin{picture}(0,0)%
\includegraphics[scale=0.3]{picture12.pstex}%
\end{picture}%
\setlength{\unitlength}{1243sp}%
\begin{picture}(11223,10853)(1318,-5942)
\put(12181,-4381){\makebox(0,0)[lb]{\smash{\SetFigFont{25}{30.0}{\familydefault}{\mddefault}{\updefault}{\color[rgb]{0,0,0}New surfaces}%
}}}
\put(12541,2744){\makebox(0,0)[lb]{\smash{\SetFigFont{25}{30.0}{\familydefault}{\mddefault}{\updefault}{\color[rgb]{0,0,0}Original surfaces}%
}}}
\put(6061,-331){\makebox(0,0)[lb]{\smash{\SetFigFont{25}{30.0}{\familydefault}{\mddefault}{\updefault}{\color[rgb]{0,0,0}time}%
}}}
\end{picture}

%% file: picture8.pstex_t
\begin{picture}(0,0)%
\includegraphics[scale=0.5]{picture8.pstex}%
\end{picture}%
\setlength{\unitlength}{2072sp}%
\begin{picture}(5872,4566)(-449,-5944)
\put(5401,-3631){\makebox(0,0)[lb]{\smash{\SetFigFont{25}{30.0}{\familydefault}{\mddefault}{\updefault}{\color[rgb]{0,0,0}quasigeodesic $\phi(\gamma_{t_1})$}%
}}}
\put(-239,-2926){\makebox(0,0)[lb]{\smash{\SetFigFont{25}{30.0}{\familydefault}{\mddefault}{\updefault}{\color[rgb]{0,0,0}geodesic arc}%
}}}
\put(-149,-5461){\makebox(0,0)[lb]{\smash{\SetFigFont{25}{30.0}{\familydefault}{\mddefault}{\updefault}{\color[rgb]{0,0,0}$\gamma_{t_2}$}%
}}}
\put(5401,-1861){\makebox(0,0)[lb]{\smash{\SetFigFont{25}{30.0}{\familydefault}{\mddefault}{\updefault}{\color[rgb]{0,0,0}$N_{K/2}(\gamma_{t_2})$}%
}}}
\put(3826,-5011){\makebox(0,0)[lb]{\smash{\SetFigFont{25}{30.0}{\familydefault}{\mddefault}{\updefault}{\color[rgb]{0,0,0}$y$}%
}}}
\put(3826,-2311){\makebox(0,0)[lb]{\smash{\SetFigFont{25}{30.0}{\familydefault}{\mddefault}{\updefault}{\color[rgb]{0,0,0}$x$}%
}}}
\end{picture}

%% file: picture4.pstex_t
\begin{picture}(0,0)%
\includegraphics[scale=0.4]{picture4.pstex}%
\end{picture}%
\setlength{\unitlength}{1658sp}%
\begin{picture}(6519,4365)(2359,-5464)
\put(5566,-4156){\makebox(0,0)[lb]{\smash{\SetFigFont{25}{30.0}{\familydefault}{\mddefault}{\updefault}{\color[rgb]{0,0,0}$\gamma_t$}%
}}}
\put(5011,-2116){\makebox(0,0)[lb]{\smash{\SetFigFont{25}{30.0}{\familydefault}{\mddefault}{\updefault}{\color[rgb]{0,0,0}$\mathcal{E}(\gamma_t)$}%
}}}
\put(4141,-5356){\makebox(0,0)[lb]{\smash{\SetFigFont{25}{30.0}{\familydefault}{\mddefault}{\updefault}{\color[rgb]{0,0,0}$\mathcal{A}(\gamma_t)$}%
}}}
\end{picture}

%% file: picture2.pstex_t
\begin{picture}(0,0)%
\includegraphics[scale=0.45]{picture2.pstex}%
\end{picture}%
\setlength{\unitlength}{1865sp}%
\begin{picture}(6082,4417)(2679,-5806)
\put(6181,-3661){\makebox(0,0)[lb]{\smash{\SetFigFont{25}{30.0}{\familydefault}{\mddefault}{\updefault}{\color[rgb]{0,0,0}$I_b$}%
}}}
\put(3931,-3631){\makebox(0,0)[lb]{\smash{\SetFigFont{25}{30.0}{\familydefault}{\mddefault}{\updefault}{\color[rgb]{0,0,0}$I_a$}%
}}}
\put(6226,-2761){\makebox(0,0)[lb]{\smash{\SetFigFont{25}{30.0}{\familydefault}{\mddefault}{\updefault}{\color[rgb]{0,0,0}$b$}%
}}}
\put(3946,-2731){\makebox(0,0)[lb]{\smash{\SetFigFont{25}{30.0}{\familydefault}{\mddefault}{\updefault}{\color[rgb]{0,0,0}$a$}%
}}}
\put(8431,-1966){\makebox(0,0)[lb]{\smash{\SetFigFont{25}{30.0}{\familydefault}{\mddefault}{\updefault}{\color[rgb]{0,0,0}original sweepout}%
}}}
\put(8761,-4186){\makebox(0,0)[lb]{\smash{\SetFigFont{25}{30.0}{\familydefault}{\mddefault}{\updefault}{\color[rgb]{0,0,0}new sweepout}%
}}}
\put(3676,-5806){\makebox(0,0)[lb]{\smash{\SetFigFont{25}{30.0}{\familydefault}{\mddefault}{\updefault}{\color[rgb]{0,0,0}Constant intervals}%
}}}
\end{picture}

%% file: picture10.pstex_t
\begin{picture}(0,0)%
\includegraphics[scale=0.27]{picture10.pstex}%
\end{picture}%
\setlength{\unitlength}{1119sp}%
\begin{picture}(20265,5521)(-899,-9835)
\put(-1690,-5461){\makebox(0,0)[lb]{\smash{\SetFigFont{25}{30.0}{\familydefault}{\mddefault}{\updefault}{\color[rgb]{0,0,0}$\hat S_t \subset \Sigma$}%
}}}
\put(-1790,-8386){\makebox(0,0)[lb]{\smash{\SetFigFont{25}{30.0}{\familydefault}{\mddefault}{\updefault}{\color[rgb]{0,0,0}$\phi(\hat S_t) \subset M$}%
}}}
\put(6776,-9736){\makebox(0,0)[lb]{\smash{\SetFigFont{25}{30.0}{\familydefault}{\mddefault}{\updefault}{\color[rgb]{0,0,0}Modified sweepout}%
}}}
\end{picture}

%% file: picture9.pstex_t
\begin{picture}(0,0)%
\includegraphics[scale=0.5]{picture9.pstex}%
\end{picture}%
\setlength{\unitlength}{2072sp}%
\begin{picture}(6558,4371)(1576,-5044)
\put(1776,-2761){\makebox(0,0)[lb]{\smash{\SetFigFont{25}{30.0}{\familydefault}{\mddefault}{\updefault}{\color[rgb]{0,0,0}Extra edges}%
}}}
\put(7976,-961){\makebox(0,0)[lb]{\smash{\SetFigFont{25}{30.0}{\familydefault}{\mddefault}{\updefault}{\color[rgb]{0,0,0}Dual graph $\tilde \Gamma$}%
}}}
\put(7876,-3661){\makebox(0,0)[lb]{\smash{\SetFigFont{25}{30.0}{\familydefault}{\mddefault}{\updefault}{\color[rgb]{0,0,0}$\tilde x_i$}%
}}}
\end{picture}

%% file: picture13.pstex_t
\begin{picture}(0,0)%
\includegraphics[scale=0.4]{picture13.pstex}%
\end{picture}%
\setlength{\unitlength}{1658sp}%
\begin{picture}(5647,5444)(6954,-3233)
\put(12601,1289){\makebox(0,0)[lb]{\smash{\SetFigFont{25}{30.0}{\familydefault}{\mddefault}{\updefault}{\color[rgb]{0,0,0}$F$}%
}}}
\put(12151,-2311){\makebox(0,0)[lb]{\smash{\SetFigFont{25}{30.0}{\familydefault}{\mddefault}{\updefault}{\color[rgb]{0,0,0}$\phi^{-1}(T) \cup D$}%
}}}
\end{picture}

%% file: picture7.pstex_t
\begin{picture}(0,0)%
\includegraphics[scale=0.4]{picture7.pstex}%
\end{picture}%
\setlength{\unitlength}{1658sp}%
\begin{picture}(9675,5966)(-899,-6680)
\put(1351,-5911){\makebox(0,0)[lb]{\smash{\SetFigFont{25}{30.0}{\familydefault}{\mddefault}{\updefault}{\color[rgb]{0,0,0}$S \times I$}%
}}}
\put(6751,-3211){\makebox(0,0)[lb]{\smash{\SetFigFont{25}{30.0}{\familydefault}{\mddefault}{\updefault}{\color[rgb]{0,0,0}$Y$}%
}}}
\put(8776,-1411){\makebox(0,0)[lb]{\smash{\SetFigFont{25}{30.0}{\familydefault}{\mddefault}{\updefault}{\color[rgb]{0,0,0}$S_n$}%
}}}
\put(4951,-6586){\makebox(0,0)[lb]{\smash{\SetFigFont{25}{30.0}{\familydefault}{\mddefault}{\updefault}{\color[rgb]{0,0,0}$C$}%
}}}
\put(8776,-3661){\makebox(0,0)[lb]{\smash{\SetFigFont{25}{30.0}{\familydefault}{\mddefault}{\updefault}{\color[rgb]{0,0,0}$T_3$}%
}}}
\put(8776,-6361){\makebox(0,0)[lb]{\smash{\SetFigFont{25}{30.0}{\familydefault}{\mddefault}{\updefault}{\color[rgb]{0,0,0}$S_1$}%
}}}
\put(8776,-5011){\makebox(0,0)[lb]{\smash{\SetFigFont{25}{30.0}{\familydefault}{\mddefault}{\updefault}{\color[rgb]{0,0,0}$S_2$}%
}}}
\put(8776,-2311){\makebox(0,0)[lb]{\smash{\SetFigFont{25}{30.0}{\familydefault}{\mddefault}{\updefault}{\color[rgb]{0,0,0}$T_{n-1}$}%
}}}
\put(3826,-3886){\makebox(0,0)[lb]{\smash{\SetFigFont{25}{30.0}{\familydefault}{\mddefault}{\updefault}{\color[rgb]{0,0,0}$\phi$}%
}}}
\put(-1099,-2986){\makebox(0,0)[lb]{\smash{\SetFigFont{25}{30.0}{\familydefault}{\mddefault}{\updefault}{\color[rgb]{0,0,0}$\phi^{-1}(Y)$}%
}}}
\end{picture}